# Apostol-Bernoulli functions, derivative polynomials and Eulerian polynomials

Khristo N. Boyadzhiev
Ohio Northern University
Departnment of Mathematics
Ada, OH 45810
k-boyadzhiev@onu.edu

**1. Introduction**

The Apostol-Bernoulli functions $\beta_n(a,\lambda)$ are defined by the exponential generating function

$$g(z,\lambda,a) = \frac{ze^{az}}{\lambda e^z - 1} = \sum_{n=0}^{\infty} \beta_n(a,\lambda) \frac{z^n}{n!} \ . \tag{1.1}$$

Thus, when $\lambda \neq 1$,

$$\beta_0(a,\lambda) = 0, \ \beta_1(a,\lambda) = \frac{1}{\lambda - 1}, \ \beta_2(a,\lambda) = \frac{2a(\lambda-1) - 2\lambda}{(\lambda-1)^2}, \ldots \text{etc} \tag{1.2}$$

These functions were introduced by Apostol [2] in order to evaluate the Lerch transcendent (also Lerch zeta function)

$$\Phi(\lambda,s,a) = \sum_{n=0}^{\infty} \frac{\lambda^n}{(n+a)^s} \tag{1.3}$$

for negative integer values of $s$ - see (6.5) below. Here $|\lambda| \leq 1$ and $a > 0$. A detailed definition of $\Phi$ and its basic properties can be found in [8].

The functions $\beta_n(a,\lambda)$ are polynomials in the first variable, $a$, and rational functions in the second variable, $\lambda$. They were studied and generalized recently in a number of papers, mainly by Luo and Srivastava, under the name Apostol-Bernoulli polynomials, [10], [12], [13].

We present here a short elementary survey with a special focus - pointing out the relation of $\beta_n(a,\lambda)$ to the classical Eulerian polynomials [7] and the derivative polynomials considered in [5] and [6].

## 2. Basic properties

We list here are some immediate properties of $\beta_n(a,\lambda)$ following from its definition (1.1); see [2], [3]. First, remind that the Bernoulli polynomials $B_n(a)$ are defined by the generating function

$$\frac{ze^{az}}{e^z - 1} = \sum_{n=0}^{\infty} B_n(a)\frac{z^n}{n!}, \qquad (2.1)$$

([1], [8], [14]). Therefore, when $\lambda = 1$, it follows from (1.1) that $\beta_n(a,1) = B_n(a)$. Because of this relationship, $\beta_n(a,\lambda)$ were called Apostol-Bernoulli polynomials by Luo [12].

For $\lambda \neq 1$ Apostol defined the functions $\beta_n(\lambda) = \beta_n(0,\lambda)$ for which we have

$$\frac{z}{\lambda e^z - 1} = \sum_{n=0}^{\infty} \beta_n(\lambda)\frac{z^n}{n!}. \qquad (2.2)$$

Note that when $n > 0$ the functions $\beta_n(\lambda)$ are not defined for $\lambda = 1$, see (2.5) below. This discontinuity follows from the discontinuity of the generating gunction $g(z,\lambda,a)$ in (1.1). When $\lambda = 1$, $\lim_{z \to 0} g(z,1,a) = 1$, and when $\lambda \neq 1$, $g(0,\lambda,a) = 0$. Thus $\lim_{\lambda \to 1} \beta_n(a,\lambda) \neq B_n(a)$.

The series (1.1) can be viewed as the product of the series (2.2) and the Taylor expansion of $e^{az}$, i.e.

$$\sum_{n=0}^{\infty} \beta_n(a,\lambda)\frac{z^n}{n!} = (\sum_{m=0}^{\infty} \beta_m(\lambda)\frac{z^m}{m!})(\sum_{k=0}^{\infty} a^k \frac{z^k}{k!}). \qquad (2.3)$$

From this equation, comparing coefficient on both sides one comes to the representation

$$\beta_n(a,\lambda) = \sum_{k=0}^{n} \binom{n}{k} \beta_k(\lambda) a^{n-k} \qquad (2.4)$$

($n = 0, 1, \ldots,$), which shows, in particular, that the functions $\beta_n(a, \lambda)$ are polynomial in the variable $a$ of order $n-1$. The functions $\beta_n(\lambda)$ are not polynomials, but rational functions,

$$\beta_0(\lambda) = 0, \; \beta_1(\lambda) = \frac{1}{\lambda - 1}, \; \beta_2(\lambda) = \frac{-2\lambda}{(\lambda - 1)^2}, \; \beta_3(\lambda) = \frac{3(\lambda^2 + \lambda)}{(\lambda - 1)^3}, \text{ etc}. \qquad (2.5)$$

These functions can be computed recursively for $n \geq 2$ from the equation

$$\beta_n(\lambda) = \lambda \sum_{k=0}^{n} \binom{n}{k} \beta_k(\lambda). \qquad (2.6)$$

Equation (2.6) follows from (2.4) and (2.16) below. Apostol computed the first six $\beta_k(\lambda)$ in [2] and also obtained the general formula

$$\beta_k(\lambda) = \frac{k}{\lambda - 1} \sum_{j=0}^{k-1} \{{}^{k-1}_j\} j! (\frac{\lambda}{1 - \lambda})^j, \qquad (2.7)$$

where $\{{}^m_k\}$ are the Stirling numbers of the second kind (see [7], [11]). Combining this with (2.4) one finds the representation (cf. [12])

$$\beta_n(a, \lambda) = \frac{1}{\lambda - 1} \sum_{k=0}^{n} \binom{n}{k} k a^{n-k} \sum_{j=0}^{k-1} \{{}^{k-1}_j\} j! (\frac{\lambda}{1 - \lambda})^j. \qquad (2.8)$$

Here are some further properties of $\beta_n(a, \lambda)$ from [2]

$$(\frac{\partial}{\partial \lambda})^p \beta_n(a, \lambda) = \frac{n!}{(n-p)!} \beta_{n-p}(a, \lambda), \; (0 \leq p \leq n), \qquad (2.9)$$

$$\beta_n(a + b, \lambda) = \sum_{k=0}^{n} \binom{n}{k} \beta_k(a, \lambda) b^{n-k}, \qquad (2.10)$$

$$\int_a^b \beta_n(t, \lambda) dt = \frac{1}{n+1} (\beta_{n+1}(b, \lambda) - \beta_{n+1}(a, \lambda)). \qquad (2.11)$$

Definition (1.1) leads also to the difference equation [2, (3.3)]

$$\lambda \beta_n(a+1,\lambda) - \beta_n(a,\lambda) = na^{n-1}, \tag{2.12}$$

from which by iteration

$$\lambda^m \beta_n(m,\lambda) - \beta_n(0,\lambda) = n \sum_{k=0}^{m-1} k^{n-1}\lambda^k, \text{ i.e.} \tag{2.13}$$

$$\lambda^m \beta_n(m,\lambda) = \beta_n(\lambda) + n \sum_{k=0}^{m-1} k^{n-1}\lambda^k, \tag{2.14}$$

for every $m = 1, 2, \ldots$. The polynomials on the right hand side in (2.13), without the factor $n$, are known as the Mirimanoff polynomials ([15, p. 504]). Taking $n = 1$ and $a = 0$ in (2.12), we see that

$$\lambda \beta_1(1,\lambda) = \beta_1(\lambda) + 1, \tag{2.15}$$

and with $a = 0$, for $n \geq 2$,

$$\lambda \beta_n(1,\lambda) = \beta_n(\lambda). \tag{2.16}$$

### 3. Geometric polynomials

In this section we show that the functions $\beta_n(\lambda)$ are closely related to the geometric polynomials

$$\omega_n(x) = \sum_{k=0}^{n} \left\{{n \atop k}\right\} k!\, x^k, \tag{3.1}$$

$\omega_0(x) = 1$, $\omega_1(x) = x$, $\omega_2(x) = 2x^2 + x$, $\omega_3(x) = 6x^3 + 6x^2 + x$, etc,

used in [6] for several applications. These polynomials describe the action of the derivative operator $(x\frac{d}{dx})^m$, $m = 0, 1, 2, \ldots$ on the function $\frac{1}{1-x}$,

$$(x\frac{d}{dx})^m \left\{\frac{1}{1-x}\right\} = \frac{1}{1-x} \omega_m\left(\frac{x}{1-x}\right), \quad m = 0, 1, 2, \ldots, \tag{3.2}$$

and when $|x| < 1$,

$$\sum_{k=0}^{\infty} k^m x^k = \frac{1}{1-x} \omega_m(\frac{x}{1-x}) . \qquad (3.3)$$

More generally, the following geometric summation formula is true.

**Proposition 3.1**. For a large class of entire functions and $|x|<1$,

$$\sum_{k=0}^{\infty} f(k) x^k = \frac{1}{1-x} \sum_{n=0}^{\infty} \frac{f^{(n)}(0)}{n!} \omega_n(\frac{x}{1-x}) , \qquad (3.4)$$

(see [6] for details). Equation (3,3) follows from (3.4) when $f(x) = x^m$.

Comparing (3.1) to (2.7) we find that

$$\beta_k(\lambda) = \frac{k}{\lambda - 1} \omega_{k-1}(\frac{\lambda}{1-\lambda}), \text{ or} \qquad (3.5)$$

$$\beta_k(\frac{z}{1+z}) = -k(z+1) \omega_{k-1}(z), \qquad (3.6)$$

for $k = 1, 2, \ldots$. This shows, in particular, that $\beta_k(\frac{z}{1+z})$ are polynomials of the variable $z$.

**Corollary 3.2.** We can write the summation formula (3.4) also in the form

$$\sum_{k=0}^{\infty} f(k) x^k = - \sum_{m=1}^{\infty} \frac{f^{(m-1)}(0)}{m!} \beta_m(x) . \qquad (3.7)$$

The Apostol-Bernoulli functions $\beta_n(a, \lambda)$ can be expressed in terms of the geometric polynomials,

$$\beta_n(a, \lambda) = \frac{1}{\lambda - 1} \sum_{k=0}^{n} \binom{n}{k} k a^{n-k} \omega_{k-1}(\frac{\lambda}{1-\lambda}), \qquad (3.8)$$

as follows from (2.8).

**4. Eulerian polynomials**

In 1755 Leonhard Euler published the book, *Insitutiones calculi differentialis...*, [9], where he systematically developed important topics in differential calculus. In Part 2, Chapter VII,

Euler evaluated sums of the form

$$1^n x + 2^n x^2 + \ldots + k^n x^k + \ldots \qquad (4.1)$$

which led him to the polynomials

$$A_0(x) = 1, A_1(x) = x, A_2(x) = x^2 + x, A_3(x) = x^3 + 4x^2 + x, \text{ etc,} \qquad (4.2)$$

presently known as Eulerian polynomials (not to be confused with the Euler polynomials, [1], [8], [14] ). In Paragraph 173 (see also 176) of his book Euler computed these polynomials up to $A_6(x)$. The Eulerian polynomials are usually written as

$$A_n(x) = \sum_{k=0}^{n} \langle {n \atop k} \rangle x^{n-k} \qquad (4.3)$$

where $\langle {n \atop k} \rangle$ are the Eulerian numbers (as defined in [11]). The infinite sum in (3.3) can be evaluated in terms of $A_n(x)$, namely,

$$\sum_{k=0}^{\infty} k^n x^k = \frac{A_n(x)}{(1-x)^{n+1}}, \ |x| < 1. \qquad (4.4)$$

The Eulerian numbers and polynomials have interesting combinatorial applications described, for instance, in [7] and [11] ( the Eulerian numbers $A(n,k)$ in [7] differ slightly from those in [11]: $A(n,k) = \langle {n \atop k-1} \rangle$ ).

The Eulerian polynomials $A_n$ are related to the geometric polynomials $\omega_n$ in a simple way, as follows from (3.3) and (4.4),

$$\omega_n(\frac{\lambda}{1-\lambda}) = \frac{A_n(\lambda)}{(1-\lambda)^n}. \qquad (4.5)$$

Essentially, these relations originate from Euler's work [9]. The polynomials $\omega_n$ appear in Part 2, Paragraph 172 (next to the Eulerian polynomials in 173), where several of them were evaluated. From (4.5) and (3.5) comes the relation

$$\beta_k(\lambda) = \frac{-k}{(1-\lambda)^k} A_{k-1}(\lambda), \qquad (4.6)$$

and therefore, (2.4) can be written also in the form

$$\beta_n(a,\lambda) = -\sum_{k=1}^{n} \binom{n}{k} \frac{k a^{n-k}}{(1-\lambda)^k} A_{k-1}(\lambda). \qquad (4.7)$$

**5. Derivative polynomials for tangents, cotangents and secants**

It is easy to see that the consecutive derivatives of $\tanh x$ are polynomials of $\tanh x$, as follows from the identity $\operatorname{sech}^2 x = 1 - \tanh^2 x$. (The same is true, of course, for $\coth x, \cot x, \tan x$). The derivative polynomials $C_m$ and $S_m$ for the hyperbolic tangent and secant are defined by the equations

$$\left(\frac{d}{dx}\right)^m \tanh x = C_m(\tanh x), \qquad (5.1)$$

$$\left(\frac{d}{dx}\right)^m \operatorname{sech} x = \operatorname{sech} x\, S_m(\tanh x), \qquad (5.2)$$

for $m = 0, 1, \ldots$, see [5]. Notice that $C_m$ is of degree $m+1$, while $S_m$ is of degree $m$. The polynomial $C_m$ is also the derivative polynomial for $\coth x$. The derivative polynomials $P_m$ and $Q_m$ for the trigonometric tangent and secant correspondingly, are defined likewise and differ only slightly from $C_m$ and $S_m$. For instance,

$$P_m(x) = -i^{m+1} C_m(ix), \qquad (5.3)$$

and the derivative polynomial for $\cot x$ is

$$(-1)^m P_m(z) = -P_m(-z) = (-i)^{m+1} C_m(ix). \qquad (5.4)$$

Thus

$$\left(\frac{d}{dx}\right)^m \cot x = (-i)^{m+1} C_m(i \cot x) \qquad (5.5)$$

All these polynomials, $C_m, S_m, P_m, Q_m$, were computed explicitly in [5]. They are closely related to the geometric polynomials $\omega_m$ and hence to the functions $\beta_n(\lambda)$ as well.

**Proposition 5.1** (see [5]) The derivative polynomials for $\tanh x$ ( $\coth x$ ) and $\operatorname{sech} x$ are given by

$$C_m(z) = (-2)^m (z+1) \omega_m(\frac{z-1}{2}), \tag{5.6}$$

$$S_m(z) = \sum_{k=0}^{m} \binom{m}{k} 2^k \omega_k(\frac{z+1}{-2}). \tag{5.7}$$

Also, in view of (3.5) we have the following.

**Corollary 5.2**

$$C_m(z) = \frac{(-2)^{m+1}}{m+1} \beta_{m+1}(\frac{z-1}{z+1}), \tag{5.8}$$

and

$$S_m(z) = \frac{1}{z-1} \sum_{k=0}^{m} \binom{m}{k} \frac{2^{k+1}}{k+1} \beta_{k+1}(\frac{z+1}{z-1}). \tag{5.9}$$

Combining (5.3) and (5.8) we find

$$P_m(z) = \frac{(-1)^m (2i)^{m+1}}{m+1} \beta_{m+1}(\frac{iz-1}{iz+1}). \tag{5.10}$$

Apostol in [4] evaluated several consecutive derivative polynomials for $\cot \pi x$ and obtained a recursion formula for them. These polynomials appear in the evaluation at the negative integers of the function

$$F(x,s) = \sum_{n=1}^{\infty} \frac{e^{2\pi i n x}}{n^s}. \tag{5.11}$$

The function $F(x,s)$ defined by (5.11) for $\operatorname{Re} s > 1$, can be continued analytically for all complex $s$. Apostol obtained the representation

$$F(x,-m) = \frac{1}{(2\pi i)^m} \frac{d^m}{dx^m}(\frac{i}{2}\cot \pi x), \tag{5.12}$$

for $m = 0, 1, 2, \ldots$ and non-integer $x$, [4, p 227]. Therefore, according to (5.4),

$$F(x, -m) = \left(\frac{i}{2}\right)^{m+1} P_m(\cot \pi x). \tag{5.13}$$

**6. Hankel integral representation and the negative values of the Lerch transcendent**

The following integral representation is true (cf. [2], [8])

$$\Phi(\lambda, s, a) = \frac{\Gamma(1-s)}{2\pi i} \int_L \frac{z^{s-1} e^{az}}{1 - \lambda e^z} dz, \tag{6.1}$$

where $L$ is the Hankel contour consisting of three parts: $L = L_- \cup L_+ \cup L_\epsilon$, with $L_-$ the "lower side" (i.e. $arg(z) = -\pi$) of the ray $(-\infty, -\epsilon)$, $\epsilon > 0$, traced left to right, and $L_+$ the "upper side" ($arg(z) = \pi$) of this ray traced right to left. Finally, $L_\epsilon = \{z = \epsilon e^{\theta i} : -\pi \le \theta \le \pi\}$ is a small circle traced counterclockwise and connecting the two sides of the ray. For the proof of this representation we start with $\text{Re } s > 1$, but then it is clear that the representation hold for all complex $s$ and thus provides a holomorphic extension of $\Phi(\lambda, s, a)$ in this variable. In particular, when $s = -m$, $m = 0, 1, \ldots$, we have

$$\Phi(\lambda, -m, a) = \frac{m!}{2\pi i} \int_L \frac{e^{az}}{1 - \lambda e^z} \frac{dz}{z^{m+1}}. \tag{6.2}$$

It is easy to see that in this case ($m$ - integer) the integrals on $L_+$ and $L_-$ cancel each-other, so the contour reduces to $L_\epsilon$ and we can write

$$\Phi(\lambda, -m, a) = \frac{m!}{2\pi i} \int_{L_\epsilon} \frac{e^{az}}{1 - \lambda e^z} \frac{dz}{z^{m+1}}. \tag{6.3}$$

At the same time, it follows from (1.1) by Cauchy's formula for Taylor coefficients that

$$\beta_n(a, \lambda) = \frac{n!}{2\pi i} \int_{L_\epsilon} \frac{z e^{az}}{1 - \lambda e^z} \frac{dz}{z^{n+1}}. \tag{6.4}$$

Now (6.4) and (6.3) produce Apostol's formula [2]

$$\Phi(\lambda, -m, a) = \frac{-\beta_{m+1}(a,\lambda)}{m+1} . \tag{6.5}$$

We can compare this to (5.13). Setting $a = 1$ and $\lambda = e^{2\pi i x}$ in the definition (1.1) of $\Phi(\lambda, s, a)$, we find that

$$\Phi(e^{2\pi i x}, 1, s) = e^{-2\pi i x} F(x, s), \tag{6.6}$$

where $F(x,s)$ is the function defined in (5.11). From here, with $s = -m$,

$$F(x,-m) = e^{2\pi i x} \Phi(e^{2\pi i x}, 1, -m), \tag{6.7}$$

and therefore, in view of (6.5),

$$F(x,-m) = \frac{-e^{2\pi i x} \beta_{m+1}(1, e^{2\pi i x})}{m+1} . \tag{6.8}$$

Using property (2.16) we write this equation in the form

$$F(x,-m) = \frac{-1}{m+1} \beta_{m+1}(e^{2\pi i x}), \tag{6.9}$$

and comparing this to (5.13) we conclude that

$$(\frac{i}{2})^{m+1} P_m(\cot \pi x) = \frac{-1}{m+1} \beta_{m+1}(e^{2\pi i x}), \tag{6.10}$$

which is, of course, equivalent to equation (5.10). It is strange that Apostol did not use the connection between equations (5.12) and (6.5), obtained in his two papers [2] and [4]. His paper [2] is not even mentioned in his later paper [4].

At the end of this section we want to poit out that equation (6.5) can be used to extend the functions $\beta_n(a,\lambda)$ to complex values of $n$. Namely, we define

$$\beta_w(a,\lambda) = -w \Phi(\lambda, 1-w, a), w \in \mathbb{C}, \tag{6.11}$$

or, in terms of Hankel integral

$$\beta_w(a,\lambda) = \frac{\Gamma(1+w)}{2\pi i} \int_L \frac{z^{-w} e^{az}}{\lambda e^z - 1} dz. \tag{6.12}$$

# 7. Integral and series representation of $\beta_n(a,\lambda)$

The Lerch transcendent $\Phi(\lambda,s,a)$ has the following Hermite representation when $0 < \lambda < 1$ and for all $s$ (equation (4) in [8, p. 28]).

$$\Phi(\lambda,s,a) = \frac{1}{2a^s} + \int_0^\infty \frac{\lambda^t}{(a+t)^s} dt + 2\int_0^\infty \sin(s\arctan\frac{t}{a} - t\log\lambda)\frac{(a^2+t^2)^{-s/2}}{e^{2\pi t}-1} dt. \quad (7.1)$$

We set here $\lambda = e^{-\alpha}$, $\alpha > 0$, $s = -n$, $n > 0$ to obtain:

$$\sum_{k=0}^\infty (k+a)^n e^{-\alpha k} = \frac{a^n}{2} + \int_0^\infty (a+t)^n e^{-\alpha t} dt + 2\int_0^\infty \sin(\alpha t - n\arctan\frac{t}{a})\frac{(a^2+t^2)^{n/2}}{e^{2\pi t}-1} dt. \quad (7.2)$$

For the left hand side, which is $\Phi(e^{-\alpha},-n,a)$, we use (6.5) to write

$$\Phi(e^{-\alpha},-n,a) = \frac{-\beta_{n+1}(a,e^{-\alpha})}{n+1}. \quad (7.3)$$

Thus $\quad (7.4)$

$$\beta_{n+1}(a,e^{-\alpha}) = -(n+1)(\frac{a^n}{2} + \int_0^\infty (a+t)^n e^{-\alpha t} dt + 2\int_0^\infty \sin(\alpha t - n\arctan\frac{t}{a})\frac{(a^2+t^2)^{n/2}}{e^{2\pi t}-1} dt)$$

**Proposition 7.1**. For $\forall \alpha \neq 0$ and $m > 1$,

$$\beta_m(e^{-\alpha}) = -\frac{m!}{\alpha^m} - 2m\int_0^\infty \cos(\alpha t - \frac{m\pi}{2})\frac{t^{m-1}}{e^{2\pi t}-1} dt. \quad (7.5)$$

This results from (7.4) by taking limits of both sides when $a \to 0^+$ and setting $n+1 = m$.

We shall obtain now a relation which explains the discontinuity in (1.1). From [8, p. 30],

$$\Phi(e^{-\alpha},-n,a) = e^{a\alpha}(\frac{n!}{\alpha^{n+1}} - \sum_{k=0}^\infty \frac{B_{n+k+1}(a)\alpha^k}{k!(n+k+1)}). \quad (7.6)$$

Thus

$$\beta_m(a, e^{-\alpha}) = -e^{\alpha a}\left(\frac{m!}{\alpha^m} - m \sum_{k=0}^{\infty} \frac{B_{m+k}(a)\,\alpha^k}{k!\,(m+k)}\right),$$

and we derive from here the following.

**Proposition 7.2.** For any $m = 1, 2\ldots,$

$$\lim_{\alpha \to 0+} \left[\beta_m(a, e^{-\alpha}) + e^{\alpha a}\frac{m!}{\alpha^m}\right] = mB_m(a), \tag{7.6}$$

$$\lim_{\alpha \to 0+} \left[\beta_m(e^{-\alpha}) + \frac{m!}{\alpha^m}\right] = mB_m, \tag{7.7}$$

where $B_m = B_m(0)$ are the Bernoulli numbers.

**References.**